\documentclass[12pt]{article}
\usepackage{latexsym}
\usepackage{amsmath}
\usepackage{amssymb}
\begin{document}

\newtheorem{theorem}{Theorem}[section]
\newtheorem{ex}[theorem]{Example}
\newtheorem{prop}[theorem]{Proposition}
\newtheorem{cor}[theorem]{Corollary}
\newtheorem{lemma}[theorem]{Lemma}
\newtheorem{define}[theorem]{Definition}
\newcommand{\qed}{\ensuremath{\hfill \Box \vspace{.2in}}}
\newcommand{\Z}{\textbf{Z}}
\newcommand{\length}{\textrm{len}}
\title{Simple Regular Skew Group Rings\thanks{The research in this paper
 will form part of the author's Ph.D. dissertation at the University of California
 at Santa Barbara.}}
\author{Kathi Crow}
\date{} \maketitle

\begin{abstract}

 Given a group $G$ acting on a ring $R$ via $\alpha:G\rightarrow
 \textup{\textrm{Aut}}(R)$, one can construct
the skew group ring $R*_{\alpha}G$.  Skew group rings have been
studied in depth, but necessary and sufficient conditions for the
simplicity of a general skew group ring are not known. In this
paper, such conditions are given for certain types of skew group
rings, with an emphasis on Von Neumann regular skew group rings.
Next the results of the first section are used to construct a
class of simple skew group rings.  In particular, we obtain more
efficient proofs of the simplicity of certain rings constructed by
H. Kambara and J. Trlifaj.
\end{abstract}

\section{Introduction}

Let $\alpha:G\rightarrow \textup{\textrm{Aut}}(R)$ be an action of
a group $G$ on a ring $R$.  Denote $\alpha(g)(r)$ by ${}^gr$ for
all $g\in G$ and $r\in R$.  If the identity is the only element of
$G$ that maps to an inner automorphism then the action is said to
be \emph{outer}. Additively, the skew group ring $R*_{\alpha}G$ is
the free left $R$-module with basis $G$. Thus elements of
$R*_{\alpha}G$ are finite sums, $\sum r_gg$ where $r_g\in R$ and
$g\in G$. Multiplication in $R*_{\alpha}G$ is given by the
multiplication in $R$ and in $G$, and by $gr={}^grg$, for $r\in R$
and $g\in G$ and then extending linearly. The \emph{support} of an
element is \mbox{supp( $\sum r_gg$ )=}$\{g\in G\mid r_g\neq 0$\}.
The \emph{length} of $\sum r_gg$ is the number of elements in
supp($\sum r_gg$), and is denoted $\length(\sum r_gg)$.

Suppose $I$ is an ideal of $R$.  Then $I$ is
\textit{$G$-invariant} if ${}^gI\subseteq I$ for every $g\in G$.
Notice that this is equivalent to ${}^gI=I$ for every $g\in G$
because ${}^gI\subseteq I$ and ${}^{g^{-1}}I\subseteq I$ together
imply that ${}^gI=I$. Define $R$ to be \textit{$G$-simple} if $R$
is nonzero and the only $G$-invariant ideals of $R$ are $0$ and
$R$.

The following proposition relating $G$-simplicity of $R$ and
simplicity of $R*_{\alpha}G$ is well known. We give the short
proof for the reader's convenience.

\begin{prop}\label{montgomery} Let $\alpha:G\rightarrow \textup{\textrm{Aut}}(R)$
be an action of a group $G$ on a ring $R$.  If $R*_{\alpha}G$ is
simple, then $R$ is $G$-simple.
\end{prop}

Proof: Suppose $I$ is a nonzero $G$-invariant ideal of $R$.  Then
$I*_{\alpha}G$ is a nonzero ideal of $R*_{\alpha}G$.  Since
$R*_{\alpha}G$ is simple, $I*_{\alpha}G=R*_{\alpha}G$. Thus
$R\subseteq I*_{\alpha}G$, and so we have $R\subseteq I$.
Therefore $R=I$ and hence $R$ is $G$-simple. \qed

Thus the $G$-simplicity of $R$ is a necessary condition for
$R*_{\alpha}G$ to be simple.  While this is not a sufficient
condition, it is in certain cases.  In section~\ref{simplicity} we
give various conditions on the ring $R$, the group $G$, and the
action $\alpha$ which will then result in the equivalence of
$G$-simplicity of $R$ and simplicity of $R*_{\alpha}G$.  In
particular, we consider $X$-outer actions and how to determine
whether an automorphism is $X$-inner.

In the following section we look at examples that arise from
letting a group act on a topological space and considering a ring
of functions on that space.  We use the results of
section~\ref{simplicity} to determine the simplicity of particular
skew group rings constructed this way.  Finally we note that rings
constructed by Kambara and Trlifaj are isomorphic to rings of this
type and then give more efficient proofs of the simplicity of
these rings.

\section{Simplicity of skew group rings}\label{simplicity}
 The goal in this section is to find cases where $G$-simplicity of
 $R$ implies simplicity of the skew group ring. First we will
 recall a lemma which is useful when showing that a given skew
 group ring is simple.

\begin{lemma} \label{g-simple} If $R$ is $G$-simple then no proper ideal of
$R*_{\alpha}G$ intersects $R$ nontrivially.
\end{lemma}

Proof: Suppose $R$ is $G$-simple and $I$ is an ideal of
$R*_{\alpha}G$ which intersects $R$ nontrivially.  Let $J=I\bigcap
R$.  If $x$ is a nonzero element of $J$, then ${}^gx=gxg^{-1}\in
J$ for every $g\in G$. So $J$ is a $G$-invariant ideal of $R$ and
hence $J=R$. Thus $1\in I$ so $R*_{\alpha}G=I$. \qed

The following proposition is probably well known, but we were
unable to locate a reference in the literature.

\begin{prop}\label{abelian} If $G$ is an abelian group with outer action $\alpha$ on $R$,
then $R*_{\alpha}G$ is simple if and only if $R$ is $G$-simple.
\end{prop}

Proof: Proposition~\ref{montgomery} gives the result that if
$R*_{\alpha}G$ is simple then $R$ is $G$-simple. Now suppose $R$
is $G$-simple and $I$ is a nonzero ideal of $R*_{\alpha}G$. Let
$x=\sum_{i=1}^nr_ig_i$ be an element of $I$ with minimal positive
length. If $n=1$, then $xg_1^{-1}$ is a nonzero element of
$R\bigcap I$. Thus the lemma yields $I=R*_{\alpha}G$.

 Now suppose $n>1$. Moreover, without
loss of generality assume each $r_i\neq 0$, the $g_i$ are
distinct, and $g_1=1$.

Since $R$ is $G$-simple and $r_1\neq 0$, the set $\{{}^gr_1\mid
g\in G\}$ generates $R$ as an ideal.  Hence
$1=\sum_j\sum_ka_{kj}{}^{g_j}r_1b_{kj}$, for some
$a_{kj},b_{kj}\in R$ and $g_j\in G$. Let
\begin{eqnarray*}
y&=&\sum_j\sum_ka_{kj}g_jxg_j^{-1}b_{kj}
 =1+\sum_j\sum_ka_{kj}(\sum_{i=2}^n{}^{g_j}r_ig_i)b_{kj}\\
&=&1+\sum_{i=2}^n(\sum_j\sum_ka_{kj}{}^{g_j}r_i{}^{g_i}b_{kj})g_i.\\
\end{eqnarray*}
Thus $\length(y)\leq\length(x)$, so $y$ is an element of $I$ of
minimal length. Hence, we may assume that
$x=1+\sum_{i=2}^nr_ig_i$.

Let $r\in R$.  Then $rx-xr$ is in $I$ and has length strictly less
than the length of $x$.  Since $x$ was chosen to have minimal
positive length, $rx-xr=0$.  Thus we have $rr_i=r_i{}^{g_i}r$ for
each $r\in R$ and $1\leq i\leq n$.

If $h\in G$, $\length(hx-xh)<\length(x)$ and hence $hx-xh=0$. So
$r_i^h=r_i$ for each $h\in G$ and each $1\leq i\leq n$. Thus for
every $i$ between 1 and $n$, $Rr_iR=Rr_i=r_iR$ is $G$-invariant.
Hence $Rr_i=R=r_iR$, so $r_i$ is a unit, and
${}^{g_i}r=r_i^{-1}rr_i$ for every $r\in R$. Thus if $n>1$,
$\alpha(g_i)$ is inner, which contradicts the action being outer.
Therefore the length of $x$ is 1, so $x\in R$ and hence
$I=R*_{\alpha}G$ as noted above. \qed

If $G$ is not abelian, the type of actions allowed must be further
restricted.  Not only do we eliminate conjugation by a unit of
$R$, but also group elements that act as conjugation by elements
of a quotient ring of $R$. An automorphism $f$ of a semiprime ring
$R$ is \emph{$X$-inner} if there exists a nonzero element $u$ in
the left Martindale quotient ring of $R$ so that $f(r)u=ur$ for
every $r\in R$.  The action $\alpha$ is said to be
\emph{$X$-outer} if the only element of $G$ that maps to an
X-inner automorphism is the identity.

Montgomery proved \cite[Lemma 3.16]{Montgomery:fixed} that if $R$
is a semiprime ring and $\alpha:G\rightarrow
\textup{\textrm{Aut}}(R)$ is an $X$-outer action, then every
nonzero ideal of $R*_{\alpha}G$ intersects $R$ nontrivially. This
is exactly what is needed to prove $R*_{\alpha}G$ is simple.
Montgomery proved an analog of the following theorem
\cite[Corollary 3.18]{Montgomery:fixed} assuming that $R$ is
simple.  The same proof works with the condition on $R$ relaxed to
$G$-simplicity.

\begin{theorem}\label{x-outer} If the action $\alpha$ of a group $G$
on a semiprime ring $R$ is $X$-outer then $R$ is $G$-simple if and
only if $R*_{\alpha}G$ is simple. \end{theorem}

Proof: Because of Proposition~\ref{montgomery}, we only need to
prove the ``only if'' statement. So assume that $R$ is $G$-simple
and let $I$ be a nonzero ideal of $R*_{\alpha}G$. Then by
Montgomery's result $I\bigcap R$ is a nontrivial ideal of $R$.  So
by Lemma~\ref{g-simple}, $I=R*_{\alpha}G$ and hence $R*_{\alpha}G$
is simple. \qed

\begin{cor} If $R$ is a commutative domain and $\alpha$ is
a faithful action of $G$, then $R$ is
$G$-simple if and only if $R*_{\alpha}G$ is simple.
\end{cor}

Proof:  If $R$ is a commutative domain, then its Martindale
quotient ring is a field.  Thus any $X$-inner automorphism must be
the identity. Because $\alpha$ is faithful, $1$ is the only
element of $G$ which $\alpha$ maps to the identity. Therefore the
action is $X$-outer and applying the theorem yields the desired
result. \qed

Because showing an action is $X$-outer will help determine if the
skew group ring is simple, it is useful to have ways of checking
whether or not an automorphism is $X$-inner.  To do this, we need
to introduce a new concept.

\begin{define} An automorphism $f$ of the ring $R$ is
called \emph{corner-inner} if there exists a nonzero idempotent
$e\in R$ so that for $e'=f^{-1}(e)$ there exist elements $u\in
eRe'$ and $v\in e'Re$ such that the following conditions hold:
\begin{itemize}
\item[(a)] $uv=e$ and $vu=e'$ \item[(b)] $f(x)=uxv$ for $x\in
e'Re'$ \item[(c)] $f^{-1}(y)=vyu$ for $y\in eRe$.
\end{itemize}
\end{define}

Notice that part $(c)$ of the definition follows from parts $(a)$
and $(b)$.  Thus when showing that an automorphism is
corner-inner, it is sufficient to only prove conditions $(a)$ and
$(b)$.

\begin{prop}
\label{x-inner} Let $R$ be a semiprime ring and $f$ an
automorphism of $R$.
\begin{description}
    \item[(a)] If $R$ is commutative then the following are equivalent:
    \begin{description}
    \item[(i)] $f$ is X-inner.
    \item[(ii)] $\textup{\textrm{ann}}_R((\textup{\textrm{id}}-f)(R))\neq 0$.
    \item[(iii)] There is a nonzero ideal $I$ in $R$ so that $f$ is the
    identity on $I$.
\end{description}
    \item[(b)] If $R$ is commutative and regular, then $f$ is
    X-inner if and only if there exists a nonzero idempotent $e\in
    R$ so that $f$ is the identity on $eR$.
    \item[(c)] If $R$ is regular and $f$ is X-inner, then $f$ is corner-inner.
\end{description}
\end{prop}
Proof: Denote the Martindale quotient ring of $R$ by $Q$.

(a) Let $R$ be a commutative ring.

((i)$\Rightarrow$(ii)) Suppose $f$ is $X$-inner. Then there is a
nonzero element $x\in Q$ so that $f(r)x=xr$ for every $r\in R$.
Because $R$ is commutative, $Q$ is also commutative.  Thus
$f(r)x=rx$, and hence $(r-f(r))x=0$ for every $r\in R$.  Therefore
$(\textup{\textrm{id}}-f)(R)x=0$.

Since $x$ is a nonzero element of the Martindale quotient ring,
there is an ideal $A$ of $R$ so that $0\neq xA\subseteq R$.
Because $xA\subseteq \textrm{ann}_R((\textrm{id}-f)(R))$and $xA$
is nonzero, $\textrm{ann}_R(\textrm{id}-f)(R)\neq 0$.

((ii)$\Rightarrow$(iii)) Suppose $x$ is a nonzero element of
$\textrm{ann}_R((\textrm{id}-f)(R))$.  Then $xr=xf(r)$ for every
$r\in R$.  Since $R$ is semiprime and commutative, $x^2\neq 0$.
Now $xf^{-1}(x)=x^2$ so $f(x)x=f(x^2)$.  Since $f(x)x=x^2$, we
have $x^2=f(x^2)$.  Let $I$ be the ideal generated by $x^2$.  Then
because \[f(x^2r)=f(x^2)f(r)=x^2f(r)=x(xr)=x^2r\] for every $r\in
R$, $f$ is the identity on $I$.

((iii)$\Rightarrow$(i))  Suppose $f$ is the identity on $I$ and
$x$ is a nonzero element of $I$.  Then for $r\in R$, we have
$f(r)x=f(rx)=xr$ and hence $f$ is $X$-inner.

\vspace{.1in} \noindent(b) Assume $R$ is a commutative regular
ring.

Clearly if $f$ is the identity on $eR$ for some nonzero idempotent
$e$, then by part (a) $f$ is $X$-inner.

Conversely, assume that $f$ is $X$-inner.  Then there is a nonzero
ideal $I$ in $R$ so that $f$ is the identity on $I$.  Because $R$
is regular, $I$ contains a nonzero idempotent $e$.  Thus $f$ is
the identity on $eR$.

\vspace{.1in} \noindent(c) Suppose $R$ is regular and $f$ is
$X$-inner. There is a nonzero element $q\in Q$ so that $f(r)q=qr$
for every $r\in R$.  Since $q$ is in $Q$, there is an essential
ideal $A$ of $R$ so that $qA$ is a nonzero right ideal of $R$. Now
$A$ and $f(A)$ are essential, so $A\bigcap f(A)\bigcap qA$ is a
nonzero right ideal of $R$ and hence contains a nonzero
idempotent, $e$. Because $e\in qA$, there exists $a\in A$ so that
$e=qa$.  Then $e=(eq)(ae)$ and both $eq$ and $ae$ are in $R$
because $eq=qf^{-1}(e)\in qA$. Define $e'=f^{-1}(e)$. Now
$eq=f(e')q=qe'$. Thus
\begin{eqnarray*}
 f(ae\cdot eq)& = &f(aeqe')=  f(aeq)e \\
 & = & f(aeq)qa  =  qaeqa=e,\\
\end{eqnarray*} so $ae\cdot eq=e'$.

Let $u=eq$ and $v=e'ae$.  Then $u=qe'$ and hence $u\in eRe'$.
Also, $v\in e'Re$, while $uv=e$ and $vu=e'$. Moreover, for $x\in
e'Re'$,
\[f(x)=ef(xe')e=ef(xe')qae=eqxe'ae=uxv.\] Therefore on $e'Re'$,
$f(x)=uxv$, so $f$ is corner-inner. \qed

Combining parts (a) and (c) of this result with Theorem
~\ref{x-outer} we get two new theorems equating $G$-simplicity of
$R$ with simplicity of the skew group ring.  Because the proof of
Theorem~\ref{x-outer} uses nontrivial results about the Martindale
quotient ring, we will give direct proofs of the theorems.

\begin{theorem}\label{identity on ideal}  Suppose $R$ is a commutative semiprime ring and
$\alpha$ is an action of a group $G$ on $R$.  Assume that $1$ is
the only element of $G$ whose image under $\alpha$ is the identity
on some nonzero ideal of $R$.  Then $R*_{\alpha}G$ is simple if
and only if $R$ is $G$-simple.\end{theorem}

Proof:  Proposition~\ref{montgomery} implies that if
$R*_{\alpha}G$ is simple then $R$ is $G$-simple.

Now suppose $R$ is $G$-simple. Let $I$ be a nonzero ideal of
$R*_{\alpha}G$ and $x$ a nonzero element of $I$ with minimal
length.  If $\length(x)=1$ then $x\in R$ and we are done by
Lemma~\ref{g-simple}. Otherwise we may assume
$x=r_1+r_2g_2+\ldots+r_ng_n$ with $n\geq 2$, each $r_i\neq 0 $,
and distinct $g_i\neq 1$.

For any $r\in R$,
\begin{eqnarray*}
rx-xr&=&(rr_1-r_1r)+(rr_2-r_2{}^{g_2}\!r)g_2+\ldots+(rr_n-r_n{}^{g_n}\!r)g_n \\
   &=&(rr_2-{}^{g_2}rr_2)g_2+\ldots+(rr_n-{}^{g_n}rr_n)g_n.   \\
  \end{eqnarray*}

Thus $\length(rx-xr)< n$ and since it is an element of $I$,
$rx-xr=0$ for any $r\in R$.  Therefore $rr_2={}^{g_2}rr_2$ for any
$r\in R$.  Hence $r_2^2={}^{g_2}r_2r_2$ and
${}^{g_2^{-1}}\!r_2r_2=r_2^2$, which implies that
${}^{g_2}r_2r_2={}^{g_2}(r_2^2)$.  Thus $r_2^2={}^{g_2}(r_2^2)$.
Now for $r\in R$,
\[{}^{g_2}(r_2^2r)=r_2^2{}^{g_2}r=r_2(r_2{}^{g_2}r)=r_2^2r.\]
Therefore $\alpha(g_2)$ is the identity on the nonzero ideal
$r_2^2R$.  Thus $g_2=1$.  But we assumed otherwise and hence have
a contradiction to $\length(x)>1$.  Therefore $x\in R$ and
$I=R*_{\alpha}G$ as above.  Thus $R*_{\alpha}G$ is simple. \qed

\begin{theorem}  Suppose $R$ is a regular ring and $\alpha$ is an
action of a group $G$ on $R$ so that $1$ is the only element of
$G$ whose image under $\alpha$ is corner-inner. Then
$R*_{\alpha}G$ is simple if and only if $R$ is $G$-simple.
\end{theorem}

Proof: If $R*_{\alpha}G$ is simple, $G$-simplicity of $R$ follows
from Proposition~\ref{montgomery}.

Suppose $R$ is $G$-simple.  Let $I$ be a nonzero proper ideal of
$R*_{\alpha}G$, and $x$ a nonzero element of $I$ with minimal
length.  If the length of $x$ is 1 then $xg^{-1}\in R$ for some
$g\in G$, and by an earlier lemma, $I=R*_{\alpha}G$. Now suppose
$\length(x)=n>1$. Then $x=\sum_{i=1}^nr_ig_i$, for some $r_i\in R$
and distinct $g_i\in G$. Since $R$ is regular, there exists $s\in
R$ so that $sr_1$ is a nonzero idempotent, say $e=sr_1$.  Then
$esxg_1^{-1}=e+\sum_{i=2}^nesr_ng_ng_1^{-1}\in I$ and
$1\leq\length(esxg_1^{-1})\leq\length(x)$.  Thus we may assume
that $x=e+\sum_{i=2}^nr_ig_i$ and $er_i=r_i$ for each $i$.

Since $\length(x)\geq 2$, we have $r_2\neq 0$ and $g_2\neq 1$.
There exists a nonzero element $y\in R$ such that $r_2yr_2=r_2$.
Now
\begin{eqnarray*}
\lefteqn{x-x\:{}^{g_2^{-1}}\!\!(yr_2)=}\\
& &(e-e\:{}^{g_2^{-1}}\!\!(yr_2))+(r_2-r_2yr_2)g_2+\dots+
  (r_n-r_n\:{}^{g_ng_2^{-1}}\!(yr_2))g_n= \\
   &&(e-e\:{}^{g_2^{-1}}\!\!(yr_2))+
  (r_3-r_3\:{}^{g_3g_2^{-1}}\!(yr_2))g_3+\dots+
  (r_n-r_n\:{}^{g_ng_2^{-1}}\!(yr_2))g_n \\
  \end{eqnarray*}
is an element of $I$ with length less than $n$, so
$x-x\:{}^{g_2^{-1}}\!\!(yr_2)=0$ and hence
$e=e\:{}^{g_2^{-1}}\!\!(yr_2)$.

Similarly $x-r_2yx=0$, so $r_2ye=e$.  Also, $xe-x=0$ and hence
$r_2\:{}^{g_2}\!e=r_2$. If $s\in eRe$, then $sx-xs=0$ and thus
$sr_2=r_2\:{}^{g_2}\!s$.

Let $e'={}^{g_2}\!e$.  Then the equations above yield $e'=e'yr_2$
and $r_2e'=r_2$.  Now define $u=r_2e'$ and $v=e'ye$. Then $u\in
eRe'$ and $v\in e'Re$.  We have
\[uv=(r_2e')(e'ye)=r_2e'ye=r_2ye=e,\] and also \[
vu=(e'ye)(r_2e')=e'yr_2e'=e'e'=e'. \]

Finally, for any $t\in e'Re'$,
\begin{eqnarray*}
utv&=&r_2tye=r_2\:{}^{g_2}\!(e\:{}^{g_2^{-1}}\!\!te)ye \\
&=&(e\:{}^{g_2^{-1}}\!\!te)r_2ye={}^{g_2^{-1}}\!\!(e'te')e\\
&=&{}^{g_2^{-1}}\!\!t=\alpha(g_2^{-1})(t). \\
\end{eqnarray*}
Thus $\alpha(g_2^{-1})$ is corner-inner.  Since $g_2^{-1}\neq 1$
we have a contradiction and therefore $\length(x)=1$ and so
$I=R*_{\alpha}G$.  Hence $R*_{\alpha}G$ is simple. \qed

 \section{Examples}

 An easy way to construct a skew group ring is to let $G$ be a group acting via
homeomorphisms on a topological space $X$. If $T$ is the ring of
locally constant functions from $X$ to a field $k$, then $G$ acts
on $T$ by composition. Thus we get a skew group ring
$T*_{\alpha}G$. Since $T$ is Von Neumann regular, as long as $G$
is a locally finite group and no subgroup of $G$ has order
divisible by char($k$), the resulting skew group ring is Von
Neumann regular \cite{ara:regular}, and hence semiprime. Since our
goal is to construct simple skew group rings, it is useful to know
what conditions are needed on a group action on a topological
space for our construction to yield a simple skew group ring.

Throughout this section, $G$ is a locally finite group of
homeomorphisms on a totally disconnected, compact, Hausdorff
topological space $X$ and $T$ is the ring of locally constant
functions from $X$ to a field $k$ of characteristic 0.  A subset
$U$ of $X$ is called \emph{clopen} if it is both open and closed.
Denote the complement of any set $V$ in $X$ by $V^c$. Clearly if
$V$ is clopen then so is $V^c$. Since $X$ is compact, functions in
$T$ have only finitely many values. If $g\in G$ and $t\in T$ then
$G$ acts on $T$ by $(\alpha(g)t)(x)={}^gt(x)=t(g^{-1}(x))$.

The ring $T$ must be $G$-simple if $T*_{\alpha}G$ is going to be
simple. Thus criteria are needed to check if $T$ is $G$-simple.
Here separate necessary conditions, and sufficient conditions for
$G$-simplicity of $T$ are given.

\begin{theorem}\label{open}
If the only open subsets of $X$ invariant under $G$ are
$\emptyset$ and $X$, then $T$ is G-simple.
\end{theorem}

Proof: Let $I$ be a nonzero $G$-invariant ideal of $T$, and $f$ a
nonzero element of $I$.  If $f$ is a unit then $I=T$.  So assume
otherwise. Because $f$ is locally constant, there exist nontrivial
pairwise disjoint clopen sets $U_0,U_1,\ldots , U_m$ so that for
each $i$, $f|_{U_i}$ is constant, $f(U_0)=0$ and
$f(U_i)=\lambda_i\neq 0$ if $i>0$, and $\bigcup_{i=0}^mU_i=X$.
Define $h:X\rightarrow k$ so that $h(U_0)=0$ and
$h(U_i)=\lambda_i^{-1}$ for $i>0$.  Then $\overline{f}=fh$ is zero
on $U_0$ and one on $\bigcup_{i=1}^m(U_i)$.  Moreover
$\overline{f}\in I$ so we may assume the image of $f$ is
$\{0,1\}$. Thus there is a nontrivial clopen subset $U$ of $X$
with $f(U)=1$ and $f(U^c)=0$. Then $\bigcup_{g\in G}g(U)$ is a
nonempty, open, $G$-invariant subset of $X$ and hence equals $X$.
Because $X$ is compact, there exist distinct
$g_1,g_2,\ldots,g_n\in G$ with $X=\bigcup_{i=1}^ng_i(U)$.

If $x$ is an element of $X$, there is a $g_j$ so that $x\in
g_j(U)$.  Thus if $\varphi =\sum_{i=1}^n{}^{g_i}f$, then for any
$x\in X$, $\varphi (x)\neq 0$ since $f(g_j^{-1}(x))=1$ for some
$j$ and $f(g_i^{-1}(x))=$ 0 or 1 for each $i$. Therefore $\varphi$
is a unit of $T$. But $I$ is $G$-invariant, so $\varphi \in I$ and
hence $I=T$. \qed

\begin{theorem}\label{clopen}
If \, $T$ is $G$-simple, then the only clopen subsets of $X$
invariant under $G$ are $\emptyset$ and $X$.
\end{theorem}

Proof: Suppose there is a nontrivial G-invariant clopen subset $A$
of X. Then define \mbox{$f:X\rightarrow k$} by $f(A)=0$ and
$f(A^c)=1$. Clearly $f\in T$.  Now $f$ is nonzero since $A^c\neq
\emptyset$, and $f$ is not a unit since $A\neq \emptyset$.  Let
$I$ be the ideal of $T$ generated by $f$.  The set $A$ is
G-invariant, so ${}^gf(x)=f(x)$ for every $x\in X$ and every $g\in
G$. Thus ${}^gf=f$, so $I$ is a nontrivial, $G$-invariant ideal of
$T$. \qed

It would be nice to have a necessary and sufficient condition for
$G$-simplicity of $T$ along the lines of the previous two
theorems.  Unfortunately this can not be done.  There are spaces
with groups acting on them such that the converse of
Theorem~\ref{clopen} does not hold.

\begin{ex}  A totally disconnected, compact, Hausdorff space
$X$ with a locally finite group $G$ acting via homeomorphisms so
that the ring $T$ of locally constant functions from $X$ to $k$ is
not $G$-simple and yet there are no nontrivial clopen subsets of
$X$ invariant under $G$.
\end{ex}

Let $X=\mathbb{N}\bigcup\{\infty\}$ be the one-point
compactification of $\mathbb{N}$ with the convention that
$\infty>n$ for every $n\in \mathbb{N}$.  Suppose $T$ is the ring
of all locally constant functions from $X$ to $k$. For each
$n\in\mathbb{N}$, define $g_n:X\rightarrow X$ by
\[ g_n(x)= \left\{
\begin{array}{ll}
x+1 & \textrm{if } x<n\\
 1& \textrm{if } x=n\\
x & \textrm{if } x>n \\
\end{array} \right. .\]

Then each $g_n$ is a homeomorphism of $X$ and $g_n$ generates a
cyclic group of order $n$. Let $G$ be the group generated by all
of the $g_n$.  Then $G$ is a locally finite group acting on $X$
via homeomorphisms. Suppose $I=\{f\in T\mid f(\infty)=0\}$. Then
$I$ is a nontrivial $G$-invariant ideal of $T$ and hence $T$ is
not $G$-simple.

Now suppose $U$ is a proper clopen subset of $X$.  If $\infty\in
U$ then there exists $m\geq 2$ so that $x\geq m$ implies $x\in U$
and yet $m-1\notin U$. There exists $n\in\mathbb{N}$ so that
$g_n^{-1}(m)=m-1$.  Thus $(g_n^{-1})(m)\notin U$ so $U$ is not
$G$-invariant. If $\infty\notin U$ then $U^c$ is a proper clopen
subset of $X$ containing $\infty$ and by the above argument, $U^c$
is not $G$-invariant.  If $U^c$ is not $G$-invariant then $U$ is
not either.  In both cases, the clopen set $U$ is not
$G$-invariant and hence no proper clopen subset of $X$ is
$G$-invariant.

Thus we have a space $X$ which is compact, Hausdorff, and totally
disconnected, with a locally finite group $G$ acting on $X$ via
homeomorphisms so that no proper clopen subset of $X$ is
$G$-invariant and yet the ring $T$ is not $G$-simple.\qed

Thus we do not have a condition that is both necessary and
sufficient for proving the $G$-simplicity of $T$.  But the
conditions we have are still very useful in proving the
$G$-simplicity of certain skew group rings.  This method provides
an easier way to approach already known examples of simple rings.

\begin{ex} \label{kambara} \end{ex}
Let $X_n=\{1,2,\ldots,2^n\}$ with the discrete topology. Then
$X_n$ can be mapped onto $X_{n-1}$ via $\phi_n(i)=\lfloor
\frac{(i+1)}{2} \rfloor$, where $\lfloor \ldots \rfloor$ denotes
the floor function. For $m<n$ define $\phi_{n,m}:X_n\rightarrow
X_m$ by $\phi_{n,m}=\phi_{m+1}\phi_{m+2}\ldots
\phi_{n-1}\phi_{n}$. Let $X$ be the inverse limit of this inverse
system with projection maps ${\pi}_n:X\rightarrow X_n$ and $T$ the
ring of all locally constant functions from $X$ to the field $k$.
Since each $X_n$ is compact, Hausdorff, and totally disconnected,
so is $X$.

Let $g_n:X_n\rightarrow X_n$ be the homeomorphism that sends $i$
to $i+1$ for $i<2^n$ and sends $2^n$ to $1$.  Notice that
$\phi_m\circ g_m^2=g_{m-1}\circ\phi_m$ for every $m\in
\mathbb{N}$. We claim that each $g_n$ extends to a map
$\overline{g_n}:X\rightarrow X$ by
\[\overline{g_n}(y)=(z_m)_{m\in \mathbb{N}} \textrm{\, where $z_m$}=\left\{
\begin{array}{ll}
g_m^{2^{m-n}}{\pi}_m(y) & \textrm{for $m\geq n$} \\
\phi_{n,m}(g_n{\pi}_n(y)) & \textrm{for $m< n$}\\
\end{array} \right. .\]
To check that $\overline{g_n}(y)\in X$, we must show that
$\phi_m(z_m)=z_{m-1}$ for every $m\in \mathbb{N}$.  If $m>n$, then
\[
\phi_m(z_m)=\phi_m(g_m^{2^{m-n}}{\pi}_m(y))=g_{m-1}^{2^{m-1-n}}\phi_m{\pi}_m(y)
=g_{m-1}^{2^{m-1-n}}{\pi}_{m-1}(y)=z_{m-1}.\]  On the other hand,
if $m<n$, then \[
\phi_m(z_m)=\phi_m(\phi_{n,m}g_n{\pi}_n(y))=\phi_{n,m-1}g_n{\pi}_n(y)=z_{m-1}.\]
Now if $m=n$, we have \[ \phi_n(z_n)=
\phi_n(g_n{\pi}_n(y))=\phi_{n,n-1}(g_n{\pi}_n(y))=z_{n-1}.\] Since
in each case $\phi_m(z_m)=z_{m-1}$, we do have
$\overline{g_n}(y)\in X$ as desired.

Moreover this definition yields a continuous function of $X$
because $\pi_m\circ \overline{g_n}= g_m^{2^{m-n}}\circ \pi_m$ or
$\pi_m\circ \overline{g_n}=\phi_{n,m}\circ g_n\circ \pi_n$, both
of which are continuous for every $n,m\in \mathbb{N}$.  Since
$\overline{g_n}$ generates a cyclic group of order $2^n$, we have
$\overline{g_n}^{-1}=\overline{g_n}^{2^n-1}$, which is continuous
and hence $\overline{g_n}$ is a homeomorphism of $X$.

Let $G_n$ be the cyclic group of order $2^n$ generated by
$\overline{g_n}$.  Now
$\pi_m\circ\overline{g_n}^2=\pi_m\circ\overline{g_{n-1}}$ for
every $m,n\in \mathbb{N}$.  So
$\overline{g_n}^2=\overline{g_{n-1}}$ for $n>1$ and hence
$G_{n-1}\subseteq G_n$ for every $n>1$. Thus if $G$ is the group
generated by all the $\overline{g_n}$, then $G$ is abelian,
locally finite, and in fact $G\cong\mathbb{Z}_{2^\infty}$.

As indicated above, this induces an action $\alpha$ on the ring
$T$.  Since $T$ is commutative, $\alpha$ is outer. Because $G$ is
an abelian, locally finite group and char$(k)=0$, we obtain a
regular skew group ring $T*_{\alpha}G$.

Suppose $U$ is a nonempty proper open subset of $X$.  Let
$y=(y_m)_{m\in \mathbb{N}}\in U$ and $z=(z_m)_{m\in \mathbb{N}}\in
U^c$.  Because $U$ is open, there is a basic open neighborhood,
$\theta\subseteq U$, of $y$ so that \[\theta=\left(\{y_1\}\times
\{y_2\}\times \ldots\times\{y_m\}\times \prod_{i>m}
X_i\right)\bigcap X.\] There exists $g_{m+1}^j\in G_m$ so that
$g_{m+1}^j(z_{m+1})=y_{m+1}$. If
$\overline{g}=\overline{g_{m+1}}^j$, then for $i\leq m$
\[\pi_i\overline{g}(z)=\phi_{m+1,i}(g_{m+1}^j\pi_{m+1}(z))=\phi_{m+1,i}(y_{m+1})=y_i.\]
Thus $\overline{g}(z)\in \theta$. Hence $U$ is not $G$-invariant
for any nonempty proper open $U\subseteq X$, so $T$ is $G$-simple
by Theorem~\ref{open}. Now $G$ is abelian, $\alpha$ is outer, and
$T$ is $G$-simple, so $T*_{\alpha}G$ is simple by
Proposition~\ref{abelian}. \qed

\begin{ex}  \end{ex}
Let $T_n$ be a direct product of $2^n$ copies of the field $k$ for
each $n$, with primitive idempotents $e_1^{(n)},\ldots,
e_{2^n}^{(n)}$.  Let $\tau_n:T_n\rightarrow T_{n+1}$ be the
$k$-algebra map so that
$\tau_n(e_i^{(n)})=e_{2i-1}^{(n+1)}+e_{2i}^{(n+1)}$ for each $i$
and let $T'$ be the direct limit of this system of $k$-algebras
with $\iota_n$ the injection from $T_n$ to $T'$. For $n\leq j$,
define $k$-algebra maps $g_{n,j}:T_j\rightarrow T_j$ by
$g_{n,j}(e_i^{(j)})=e_{i+2^{(j-n)}}^{(j)}$ for $i\leq
2^j-2^{(j-n)}$ and
$g_{n,j}(e_{i}^{(j)})=e_{i-2^j+2^{(j-n)}}^{(j)}$ for
$i>2^j-2^{(j-n)}$. Then each $g_{n,j}$ is an automorphism of
$T_j$. Moreover $\tau_ig_{n,i}=g_{n,i+1}\tau_i$ for every $i\in
\mathbb{N}$. Thus for each $n\in\mathbb{N}$ there exists an
automorphism $g_n$ of $T'$ so that for every $i\geq n$,
$g_n\circ\iota_i=\iota_i\circ{g_{n,i}}$. If $G'$ is the group
generated by all of the $g_{n}$ then
$G'\cong\mathbb{Z}_{2^{\infty}}$ because each $g_{n}$ generates a
cyclic group of order $2^n$ and for each $n$, $g_{n+1}^2=g_n.$
Since $G'$ is a subgroup of $\textup{\textrm{Aut}}(T')$, we have a
group action $\beta:G'\rightarrow \textup{\textrm{Aut}}(T')$. If
$\overline{T'}$ is the maximal right quotient ring of $T'$, then
every automorphism of $T'$ extends uniquely to an automorphism of
$\overline{T'}$. Thus we have a group action
$\overline{\beta}:G'\rightarrow\textup{\textrm{Aut}}(\overline{T'})$.
Since both $T'$ and $\overline{T'}$ are commutative, $\beta$ and
$\overline{{\beta}}$ are outer.

Now $G$ is isomorphic to the group $G'$ from the previous example.
If we identify $G$ and $G'$, then $T'$ is isomorphic to the ring
$T$ from that example via a $G$-equivariant isomorphism. Therefore
$T'$ is $G'$-simple. If $I$ is a nonzero $G'$-invariant ideal of
$\overline{T'}$ then $I\bigcap T'$ is a nonzero $G'$-invariant
ideal of $T'$.  But since $T'$ is $G'$-simple, we have $I\bigcap
T'=T'$ and so $I=\overline{T'}$. Thus $\overline{T'}$ is
$G'$-simple, so $\overline{T'}*_{\overline{\beta}}G'$ is simple by
Proposition~\ref{abelian}.

 Moreover
$\overline{T'}*_{\overline{\beta}}G'$ is isomorphic to a ring that
Kambara~\cite[page 112]{Kambara:existence} constructed.  He showed
that the right maximal quotient ring of this ring $R$ is directly
finite and only one-sided self-injective.  Our approach provides
an alternate proof of the simplicity of $R$, which is a key
ingredient in Kambara's construction.\qed

Another nice class of examples arises if $X$ is a topological
group and $G$ is a dense subgroup of $X$ acting by left
multiplication.  Given suitable conditions on the topology of $X$,
we can show that $T*_{\alpha}G$ will be simple.

\begin{prop}\label{topological group}
Suppose $X$ is a compact, Hausdorff, totally disconnected
topological group, and $G$ is a dense subgroup of $X$.  Let $G$
act on $X$ by left multiplication.  Then the skew group ring
$T*_{\alpha}G$ is simple.
\end{prop}

Proof: Suppose $V$ is a nontrivial proper open subset of $X$ and
$x\in X\backslash V$.  Then right multiplication by $x$ is a
homeomorphism of $X$ and thus $Gx$ is a dense subset of $X$. Since
$V$ is open, there is an element $g\in G$ so that $gx\in V$.  Thus
$V$ is not invariant under $G$ so by Theorem~\ref{open}, $T$ is
$G$-simple.

Now suppose $I$ is a nonzero ideal of $T$ and $1\neq g\in G$.  Let
$t$ be a nonzero element of $I$ and let $U$ be a clopen subset of
$X$ so that $t(U)\neq 0$. Since $G$ is dense, there is an element
$h\in U\bigcap G$. Then $gh\neq h$ so there are disjoint clopen
sets $W$ and $\theta$ so that $h\in W\subseteq U$ and
$gh\in\theta$.  Let $t_0$ be a locally constant function so that
$t_0(W)=1$ and $t_0(\theta)=0$.  Then not only is $t_0t$ an
element of $I$, but $t_0t(gh)=0$ and $t_0t(h)\neq 0$.  Thus
${}^{g^{-1}}(t_0t)\neq t_0t$, so $\alpha(g^{-1})$ is not the
identity on $I$.  Because $T$ is $G$-simple and $1$ is the only
element of $G$ whose image under $\alpha$ is the identity on some
nonzero ideal of $T$, by Theorem~\ref{identity on ideal}
$T*_{\alpha}G$ is simple.\qed

\begin{ex} \end{ex}
Let $G$ be a residually finite group and $\{H_i\}_{i\in I}$ a
directed system of normal subgroups of $G$ with finite index so
that $i\leq j$ implies $H_j\subseteq H_i$ and $\bigcap H_i=1$.
Then the groups $G/H_i$ together with projection maps
$\pi_{i,j}:G/H_{i}\rightarrow G/H_{j}$ for each $i\geq j$ form an
inverse system and we let $X$ be the inverse limit. If we endow
each $G/H_i$ with the discrete topology then since each $G/H_i$ is
finite and hence compact, $X$ is compact. Moreover $X$ is a
totally disconnected, Hausdorff topological group.

The image of the natural embedding of $G$ into $X$ is dense in
$X$.  Thus if we identify $G$ with this image we obtain a simple
skew group ring $T*_{\alpha}G$ as in Proposition~\ref{topological
group}.

 The ring $D_G$ derived from $KG$ described in \cite[page
2241]{Trlifaj:rings} is isomorphic to the ring $T*_{\alpha}G$.
Thus we have a different proof of the simplicity of $D_G$
\cite[Lemma 1.2]{Trlifaj:rings}. Trlifaj showed that if $G$ is a
locally finite direct product of countably many finite groups,
then this ring is a counterexample to the conjecture ``If $R$ is a
regular ring which is not semisimple then for each simple left
$R$-module $J$ there exists a nonzero right $R$-module $M$ so that
$M\otimes_RJ=0$'' \cite[Theorem 3.8]{Trlifaj:rings}. \qed

 \nocite{*}
\bibliography{mybib}

\begin{thebibliography}{1}

\bibitem{ara:regular}
Ricardo Alfaro, Pere Ara, and Angel Del~Rio.
\newblock Regular skew group rings.
\newblock {\em Journal of the Australian Mathematical Society}, 58:167--182,
  1995.

\bibitem{Kambara:existence}
Hikoji Kambara.
\newblock On existence of directly finite, only one-sided self-injective
  regular rings.
\newblock {\em Journal of Algebra}, 144(1):110--116, 1991.

\bibitem{Montgomery:fixed}
Susan Montgomery.
\newblock {\em Fixed Rings of Finite Automorphism Groups of Associative Rings},
  volume 818 of {\em Lecture Notes in Mathematics}.
\newblock Springer-Verlag, 1980.

\bibitem{Trlifaj:rings}
Jan Trlifaj.
\newblock Rings derived from group rings.
\newblock {\em Communications in Algebra}, 20(8):2239--2252, 1992.

\end{thebibliography}
\bibliographystyle{plain}

\vspace{.2in} Department of Mathematics, University of California,
Santa Barbara, CA 93106.

 \emph{E-mail address:} crow@math.ucsb.edu
\end{document}